\documentclass[letterpaper, 10pt, conference]{ieeeconf}
\IEEEoverridecommandlockouts       
\usepackage{amssymb,amsmath,amsfonts}
\usepackage{algorithm,algorithmic}
\usepackage{cite}
\usepackage{float}
\usepackage{epsfig}
\usepackage{graphicx}
\usepackage{graphics}
\usepackage{subfigure}
\usepackage{url}
\usepackage{xspace}
\usepackage[usenames,dvipsnames]{color}
\usepackage{soul}
\usepackage{mathtools}

\floatstyle{ruled}
\newfloat{model}{H}{mod}
\floatname{model}{\footnotesize Model}
\newfloat{notatio}{H}{not}
\floatname{notatio}{\footnotesize Notation}

\newenvironment{varalgorithm}[1]
  {\algorithm\renewcommand{\thealgorithm}{#1}}
  {\endalgorithm}

\newenvironment{list4}{
	\begin{list}{$\bullet$}{%
			\setlength{\itemsep}{0.05cm}
			\setlength{\labelsep}{0.2cm}
			\setlength{\labelwidth}{0.3cm}
			\setlength{\parsep}{0in} 
			\setlength{\parskip}{0in}
			\setlength{\topsep}{0in} 
			\setlength{\partopsep}{0in}
			\setlength{\leftmargin}{0.16in}}}
	{\end{list}}

\newenvironment{list4a}{
	\begin{list}{$\bullet$}{%
			\setlength{\itemsep}{0.05cm}
			\setlength{\labelsep}{0.2cm}
			\setlength{\labelwidth}{0.3cm}
			\setlength{\parsep}{0in} 
			\setlength{\parskip}{0in}
			\setlength{\topsep}{0in} 
			\setlength{\partopsep}{0in}
			\setlength{\leftmargin}{0.16in}}}
	{\end{list}}

\usepackage{url,changebar,bm,xspace,dsfont}
\let\mathbb=\mathds 
\def\day{\mathrm{day}} 

\newcommand{\T}{^{\mbox{\tiny T}}}

\newtheorem{theorem}{Theorem}

\newtheorem{assum}{Assumption}

\newtheorem{remark}{Remark}

\usepackage{amsmath,amsfonts,amssymb,amscd}
\usepackage{capt-of}
\usepackage{color}
\usepackage{cite}
\usepackage{verbatim}
\usepackage{hyperref}

\hypersetup{
    colorlinks = true,
    linkcolor = [rgb]{0,0,1},
    anchorcolor = [rgb]{0,0,1},
    citecolor = [rgb]{0.9,0.5,0},
    filecolor = [rgb]{0,0,1},
    urlcolor = [rgb]{0.9,0.5,0},
    bookmarksopen = true,
    bookmarksnumbered = true,
    breaklinks = true,
    linktocpage,
    colorlinks = true,
    linkcolor = [rgb]{0.2,0.6,0.2},
    urlcolor  = [rgb]{0,0,1},
    citecolor = [rgb]{0.9,0.5,0},
    anchorcolor = [rgb]{0.2,0.6,0.2},
}

\DeclareMathOperator*{\argmin}{arg\,min}

\begin{document}

\title{\LARGE \bf Asynchronous Distributed Optimization via ADMM \\ with Efficient Communication}

\author{Apostolos~I.~Rikos, Wei~Jiang, Themistoklis~Charalambous, and Karl~H.~Johansson
\thanks{Apostolos~I.~Rikos is with the Department of Electrical and Computer Engineering, Division of Systems Engineering, Boston University, Boston, MA 02215, US. E-mail: {\tt arikos@bu.edu}.}
\thanks{Wei Jiang resides in Hong Kong, China. Email: {\tt wjiang.lab@gmail.com}.}
\thanks{T.~Charalambous is with the Department of Electrical and Computer Engineering, School of Engineering, University of Cyprus, 1678 Nicosia, Cyprus.  
He is also with the Department of Electrical Engineering and Automation, School of Electrical Engineering, Aalto University, Espoo, Finland. Email: {\tt charalambous.themistoklis@ucy.ac.cy}.
}
\thanks{K.~H.~Johansson is with the Division of Decision and Control Systems, KTH Royal Institute of Technology, SE-100 44 Stockholm, Sweden. He is also affiliated with Digital Futures. E-mail: {\tt kallej@kth.se}.}
\thanks{Part of this work was supported by the Knut and Alice Wallenberg Foundation, the Swedish Research Council, and the Swedish Foundation for Strategic Research. 
The work of T. Charalambous was partly supported by the European Research Council (ERC) Consolidator Grant MINERVA (Grant agreement No. 101044629).
}
}

\maketitle
\pagestyle{empty}

%
%
%
%
\begin{abstract}
In this paper, we focus on an asynchronous distributed optimization problem. 
In our problem, each node is endowed with a convex local cost function, and is able to communicate with its neighbors over a directed communication network.  
Furthermore, we assume that the communication channels between nodes have limited bandwidth, and each node suffers from processing delays. 
We present a distributed algorithm which combines the Alternating Direction Method of Multipliers (ADMM) strategy with a finite time quantized averaging algorithm. 
In our proposed algorithm, nodes exchange quantized valued messages and operate in an asynchronous fashion. 
More specifically, during every iteration of our algorithm each node (i) solves a local convex optimization problem (for the one of its primal variables), and (ii) utilizes a finite-time quantized averaging algorithm to obtain the value of the second primal variable (since the cost function for the second primal variable is not decomposable). 
We show that our algorithm converges to the optimal solution at a rate of $O(1/k)$ (where $k$ is the number of time  steps) for the case where the local cost function of every node is convex and not-necessarily differentiable. 
Finally, we demonstrate the operational advantages of our algorithm against other algorithms from the literature. 
\end{abstract}


%
%
%
%
\section{Introduction}\label{sec:intro}

The problem of distributed optimization has received extensive attention in recent years. 
Due to the rise of large-scale machine learning \cite{2020:Nedich}, control \cite{SEYBOTH:2013}, and other data-driven applications \cite{2018:Stich_Jaggi}, there is a growing need to solve optimization problems that involve massive amounts of data. 
Solving these problems in a centralized way is proven to be infeasible since it is difficult or impossible to store and process large amounts of data on a single node. 


{Distributed optimization is a method that distributes data across multiple nodes. 
Each node performs computations on its stored data and collaborates with others to solve the optimization problem collectively. 
This approach optimizes a global objective function by combining each node's local objective function and coordinating with the network. 
The advantage is reducing computational and storage requirements for individual nodes. 
However, frequent communication with neighboring nodes is necessary to update optimization variables. 
This can become a bottleneck with increasing nodes or data. 
To address this issue, recent attention from the scientific community focuses on developing optimization algorithms with efficient communication. 
This leads to enhancements on scalability and operational efficiency, while mitigating issues like network congestion, latency, and bandwidth limitations.}

\noindent
{\textbf{Existing Literature.}
Most works in the literature assume that nodes can process and exchange real values. 
This may result in communication overhead, especially for algorithms requiring frequent and complex communication (see, e.g., \cite{2017:Makhdoumi_Ozdaglar, 2009:Nedic_Optim, 2021:Tiancheng_Uribe, 2016:Ling_Yongmei, 2021:Wei_Themis, 2021:Bastianello_Todescato, 2022:Khatana_Salapaka}). 
In practical applications, nodes must exchange quantized messages to efficiently utilize network resources like energy and processing power. 
For this reason, recent research focuses on communication-efficient algorithms (e.g., \cite{2013:Xavier_Markus, 2016:Shengyu_Biao, 2016:Ling_Yongmei, 2017:Tsai_Chang, 2021:Tiancheng_Uribe, 2022:Yaohua_Ling, rikos2023distributed}), but they often assume perfectly synchronized nodes or bidirectional communication, limiting their applicability. 
Addressing communication overhead remains a key challenge, necessitating the development of communication-efficient algorithms that can operate over directed networks asynchronously. 
Therefore, continued research in this area is crucial to overcoming this bottleneck and enhancing the performance of distributed optimization methods.}


\noindent
\textbf{Main Contributions.} 
Existing algorithms in the literature often assume that nodes can exchange precise values of their optimization variables and operate synchronously. 
However, transmitting exact values (often irrational numbers) necessitates an infinite number of bits and becomes infeasible. 
Moreover, synchronizing nodes within a distributed network involves costly protocols, time-consuming to execute. \
In this paper, we present a distributed optimization algorithm, which aims to address these challenges. 
More specifically, we make the following contributions. 
\\ \noindent \textbf{A.} We present a distributed optimization algorithm that leverages the advantages of the ADMM optimization strategy and operates over a directed communication graph. 
Our algorithm allows nodes to operate in an asynchronous fashion, and enables efficient communication as nodes communicate with quantized messages; see Algorithm~\ref{alg1}. 
\\ \noindent \textbf{B.} We prove that our algorithm converges to the optimal solution at a rate of $O(1/k)$ even for non-differentiable and convex local cost functions (as it is the case for similar algorithms with real-valued states). 
This rate is justified in our simulations in which our algorithm exhibits comparable performance with real-valued communication algorithms while guaranteeing efficient (quantized) communication among nodes; see Section~\ref{sec:results}.  
Furthermore, we show that the optimal solution is calculated within an error bound that depends on the quantization level; see Theorem~\ref{converge_Alg1}.

\section{NOTATION AND PRELIMINARIES}\label{sec:preliminaries}
\textbf{Notation.}
The sets of real, rational, integer and natural numbers are denoted by $ \mathbb{R}, \mathbb{Q}, \mathbb{Z}$ and $\mathbb{N}$, respectively. 
The symbol $\mathbb{Z}_{\geq 0}$ ($\mathbb{Z}_{> 0}$) denotes the set of nonnegative (positive) integer numbers. 
The symbol $\mathbb{R}_{\geq 0}$ ($\mathbb{R}_{> 0}$) denotes the set of nonnegative (positive) real numbers. 
The symbol $\mathbb{R}^n_{\geq 0}$ denotes the nonnegative orthant of the $n$-dimensional real space $\mathbb{R}^n$. 
Matrices are denoted with capital letters (e.g., $A$), and vectors with small letters (e.g., $x$). 
The transpose of matrix $A$ and vector $x$ are denoted as $A^\top$, $x^\top$, respectively. 
For any real number $a \in \mathbb{R}$, the floor $\lfloor a \rfloor$ denotes the greatest integer less than or equal to $a$ while the ceiling $\lceil a \rceil$ denotes the least integer greater than or equal to $a$. 
For any matrix $A \in \mathbb{R}^{n \times n}$, the $a_{ij}$ denotes the entry in row $i$ and column $j$. 
By $\mathbb{1}$ and $\mathbb{I}$ we denote the all-ones vector and the identity matrix of appropriate dimensions, respectively.
By $\| \cdot \|$, we denote the Euclidean norm of a vector.


\textbf{Graph Theory.} 
The communication network is captured by a directed graph (digraph) defined as $\mathcal{G} = (\mathcal{V}, \mathcal{E})$. 
This digraph consists of $n$ ($n \geq 2$) agents communicating only with their immediate neighbors, and is static (i.e., it does not change over time). 
In $\mathcal{G}$, the set of nodes is denoted as $\mathcal{V} =  \{ v_1, v_2, ..., v_n \}$, and the set of edges as $\mathcal{E} \subseteq \mathcal{V} \times \mathcal{V} \cup \{ (v_i, v_i) \ | \ v_i \in \mathcal{V} \}$ (note that each agent has also a virtual self-edge). 
The cardinality of the sets of nodes, edges are denoted as $| \mathcal{V} |  = N$, $| \mathcal{E} | = m$, respectively. 
A directed edge from node $v_i$ to node $v_l$ is denoted by $(v_l, v_i) \in \mathcal{E}$, and captures the fact that node $v_l$ can receive information from node $v_i$ at time step $k$ (but not the other way around). 
The subset of nodes that can directly transmit information to node $v_i$ is called the set of in-neighbors of $v_i$ and is represented by $\mathcal{N}_i^- = \{ v_j \in \mathcal{V} \; | \; (v_i, v_j)\in \mathcal{E}\}$. 
The subset of nodes that can directly receive information from node $v_i$ is called the set of out-neighbors of $v_i$ and is represented by $\mathcal{N}_i^+ = \{ v_l \in \mathcal{V} \; | \; (v_l, v_i)\in \mathcal{E}\}$. 
The \textit{in-degree}, and \textit{out-degree} of $v_j$ are denoted by $\mathcal{D}_i^- = | \mathcal{N}_i^- |$, $\mathcal{D}_i^+ = | \mathcal{N}_i^+ |$, respectively. 
The diameter $D$ of a digraph is the longest shortest path between any two nodes $v_l, v_i \in \mathcal{V}$. 
A directed \textit{path} from $v_i$ to $v_l$ of length $t$ exists if we can find a sequence of agents $i \equiv l_0,l_1, \dots, l_t \equiv l$ such that $(l_{\tau+1},l_{\tau}) \in \mathcal{E}$ for $ \tau = 0, 1, \dots , t-1$. 
A digraph is \textit{strongly connected} if there exists a directed path from every node $v_i$ to every node $v_l$, for every $v_i, v_l \in \mathcal{V}$. 


\textbf{ADMM Algorithm.} 
The standard ADMM algorithm \cite{2014:Wei_Wotao} is designed to solve the following problem: 
\begin{align}\label{admm_prob_min}
\min_{x \in \mathbb{R}^p, z \in \mathbb{R}^m}~ & f(x) + g(x),  \\
\text{s.t.}~ & A x + B z  =  c, \nonumber
\end{align}
where $A \in \mathbb{R}^{q \times p}$, $B \in \mathbb{R}^{q \times m}$ and $c \in \mathbb{R}^{q}$ (for $q, p, m \in \mathbb{N}$). 
In order to solve \eqref{admm_prob_min}, the augmented Lagrangian is: 
\begin{align}\label{admm_Lagrangian}
L_{\rho} (x, z, \lambda) = f(x) + g(x) & + \lambda (A x + B z  -  c) \nonumber
\\ & + \frac{\rho}{2} \| A x + B z  -  c \|^2, 
\end{align}
where $\lambda \in \mathbb{R} $ is the Lagrange multiplier, and $\rho \in \mathbb{R} $ is the positive penalty parameter. 
The primary variables $x$, $z$ and the Lagrangian multiplier $\lambda$ are initialized as $[ x, z, \lambda ]^\top = [ x^{[0]}, z^{[0]}, \lambda^{[0]} ]^\top$. 
Then, during every ADMM time step, the $x$, $z$ and $\lambda$ are updated as: 
\begin{align}
x^{[k+1]} = & \argmin_{x \in \mathbb{R}^p} L_{\rho} (x, z^{[k]}, \lambda^{[k]}), \label{admm_iterations_x} \\
z^{[k+1]} = & \argmin_{z \in \mathbb{R}^m} L_{\rho} (x^{[k+1]}, z, \lambda^{[k]}), \label{admm_iterations_z} \\ 
\lambda^{[k+1]} = & \lambda^{[k]} + \rho (A x^{[k+1]} + B z^{[k+1]}  -  c) , \label{admm_iterations_lambda} 
\end{align}
where $\rho$ in \eqref{admm_iterations_lambda} is the penalty parameter from \eqref{admm_Lagrangian}.



\textbf{Asymmetric Quantizers.}
Quantization is a strategy that lessens the number of bits needed to represent information. 
This reduces the required communication bandwidth and increases power and computation efficiency. 
Quantization is mainly used to describe communication constraints and imperfect information exchanges between nodes \cite{2019:Wei_Johansson}. 
The three main types of quantizers are (i) asymmetric, (ii) uniform, and (iii) logarithmic. 
In this paper, we rely on asymmetric quantizers to reduce the required communication bandwidth. 
{Note that the results of this paper are transferable to other quantizer types (e.g., logarithmic or uniform).} 
{Asymmetric quantizers are defined as 
\begin{equation}\label{asy_quant_defn}
    q_{\Delta}^a(\xi) = \Bigl \lfloor \frac{\xi}{\Delta} \Bigr \rfloor, 
\end{equation}
where $\Delta \in \mathbb{Q}$ is the quantization level, $\xi \in \mathbb{R}$ is the value to be quantized,} and $q_{\Delta}^a(\xi) \in \mathbb{Q}$ is the quantized version of $\xi$ with quantization level $\Delta$ {(note that the superscript ``$a$'' indicates that the quantizer is asymmetric).} 


%
%
%
%
\section{Problem Formulation}\label{sec:probForm}

\textbf{Problem Statement.} Let us consider a distributed network modeled as a digraph $\mathcal{G} = (\mathcal{V}, \mathcal{E})$ with $n  = | \mathcal{V} |$ nodes. 
In our network $\mathcal{G}$, we assume that the communication channels among nodes have limited bandwidth. 
Each node $v_i$ is endowed with a scalar local cost function $f_i(x): \mathbb{R}^p \mapsto \mathbb{R}$ only known to node $v_i$. 
In this paper we aim to develop a distributed algorithm which allows nodes to cooperatively solve the following optimization problem
\begin{align}\label{admm_prob_min_probForm}
\min_{x \in \mathbb{R}^p}~ & \sum_{i=1}^n f_i(x) , 
\end{align}
where $x \in \mathbb{R}^p$ is the global optimization variable (or common decision variable). 
We will solve \eqref{admm_prob_min_probForm} via the distributed ADMM strategy. 
Furthermore, in our solution we guarantee efficient communication between nodes (due to communication channels of limited bandwidth in the network). 

\vspace{.5cm}

\textbf{Modification of the Optimization Problem.} 
In order to solve \eqref{admm_prob_min_probForm} via the ADMM  and guarantee efficient communication between nodes, we introduce (i) the variable $x_i$ for every node $v_i$, (ii) the constraint $| x_i - x_j | \le \epsilon$ for every $v_i, v_j \in \mathcal{V}$ (where $\epsilon \in \mathbb{R}$ is an error tolerance which is predefined), and (iii) the constraint that nodes communicate with quantized values. 
The second constraint is introduced to allow an asynchronous implementation of the distributed ADMM strategy, and the third constraint to guarantee efficient communication between nodes. 
Considering the aforementioned constraints (i), (ii) and (iii),  \eqref{admm_prob_min_probForm} becomes: 
\begin{align}
\min_{x_i}~ & \sum_{i=1}^n f_i(x_i) , i=1, ..., n\label{admm_prob_min_probForm_1}
\\ \text{s.t.}~ & | x_i - x_j | \le \epsilon, \forall v_i, v_j \in \mathcal{V}, \label{admm_prob_min_probForm_2}
\\ & \text{nodes communicate with quantized values.} \label{constr_quant} 
\end{align}
Let us now define a closed nonempty convex set $\mathcal{C}$  as 
\begin{equation}\label{setC}
	\mathcal{C} = \left\{\begin{bmatrix} x_1\T & x_2\T & \ldots & x_n\T\end{bmatrix}\T \in \mathbb{R}^{np} \, :\,  \| x_i - x_j\| \le \epsilon \right\} .
\end{equation}
Furthermore, denote $ X \coloneqq \begin{bmatrix} x_1\T & x_2\T & \ldots & x_n\T\end{bmatrix}\T $ and its copy variable $ z \in \mathbb{R}^{np} $. 
This means that $\eqref{admm_prob_min_probForm_2}$ and \eqref{setC} become
\begin{equation}\label{problem_reformulated2} 
	X = z, \, \forall z \in \mathcal{C}. 
\end{equation}
Now let us define the indicator function $g(z)$ of set $\mathcal{C}$ as 
\begin{equation}\label{indicator_function_g}
	g(z) =
	\left\{ 
	\begin{array}{l}
	\begin{aligned}
	&0,\quad \text{if} \, z \in \mathcal{C},\\
	&\infty, \, \ \text{otherwise}.
	\end{aligned}
	\end{array}
	\right. 
\end{equation}
Incorporating \eqref{problem_reformulated2} and \eqref{indicator_function_g} into \eqref{admm_prob_min_probForm_1}, we have that \eqref{admm_prob_min_probForm_1} becomes 
\begin{equation}\label{objective_function_no_quant}
	\begin{aligned}
	\min_{z, x_i} \, & \left\{ \sum_{i=1}^{n} f_i(x_i) + g(z) \right\}, i=1,\ldots, n\\ \text{s.t.} \, & X - z = 0, \, \forall z \in \mathcal{C},  
 \\ & \text{nodes communicate with quantized values.} 
 \end{aligned}
	\end{equation}
As a result, in this paper, we aim to design a distributed algorithm that solves \eqref{objective_function_no_quant} via the distributed ADMM strategy.

%
%
%
%
\section{Preliminaries on Distributed Coordination}\label{sec:prelim_distrCoord}

We now present a definition of asynchrony (borrowed from \cite{2021:Wei_Themis}) that defines the operation of nodes in the network. 
Furthermore, we present a distributed coordination algorithm that operates with quantized values and is necessary for our subsequent development. 

\subsection{Definition of Asynchronous Operation}\label{Prel_asynch_oper}

\begin{figure}[ht]
        \begin{center}
		\includegraphics[width=0.9\columnwidth]{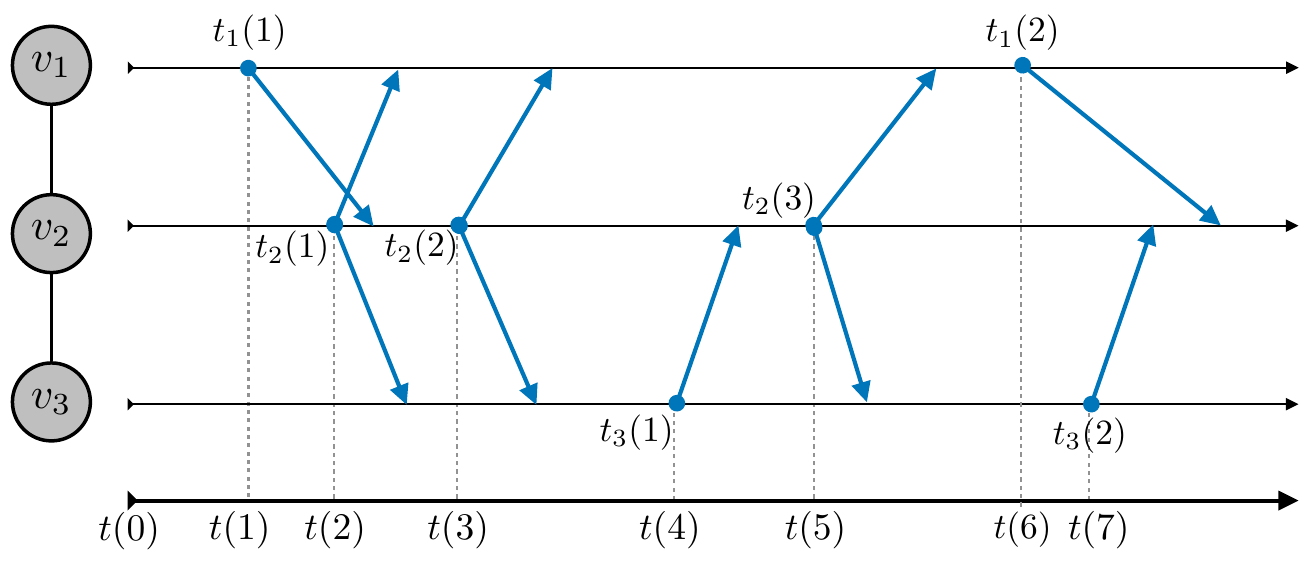}
        \end{center}
		\caption{Example of how processing and transmission delays affect the operation of nodes $v_1$, $v_2$, $v_3$. 
        Blue dots indicate the iterations and blue arrows indicate the transmissions. 
        Transmissions occur at time steps $t_i(\eta)$, and $t_i(\eta+1) - t_i(\eta)$ is the processing delay, where $i \in \{1, 2, 3\}$, $\eta \in \mathbb{Z}_{\geq 0}$. 
        The time difference from the blue dot to the blue arrow is the transmission delay \cite{2021:Wei_Themis}.} 
		\label{fig:async}
\end{figure}

During their optimization operation, nodes aim to coordinate in an asynchronous fashion. 
Specifically, let us assume that the iterations for the optimization operation start at time step $t(0)\in \mathbb{R}_{+}$. 
Furthermore, we assume that one (or more) nodes transmit values to their out-neighbors at a set of time instances $\mathcal{T}=\{t(1), t(2),t(3),\ldots\}$. 
During the nodes' asynchronous operation, a message that is received at time step $t(\eta_1)$ from node $v_i$, is processed at time step $t(\eta_2)$ where $\eta_2 > \eta_1$. 
This means that the message received at time step $t(\eta_1)$ suffers from a processing delay of $t(\eta_2)-t(\eta_1)$ time steps. 
An example of how processing delays affecting transmissions is shown in Fig.~\ref{fig:async} (that is borrowed from \cite{2021:Wei_Themis}). 

Note here that the nodes states at time step $t(\eta)$ are indexed by $\eta$. 
This means that the state of node $v_i$ at time step $t(\eta)$ is denoted as $x_i^{\eta}  \in \mathbb{R}^p$. 

We now present the following assumption which is necessary for the asynchronous operation of every node. 

\begin{assum}\label{assump_time_index_bound}
The number of time steps required for a node $v_i$ to process the information received from its in-neighbors is upper bounded by $\mathcal{B} \in \mathbb{N}$. 
Furthermore, the actual time (in seconds) required for a node $v_i$ to process the information received from its in-neighbors is upper bounded by $T \in \mathbb{R}_{\geq 0}$. 
\end{assum}

{Assumption~\ref{assump_time_index_bound} states that there exists a finite number of steps $\mathcal{B}$ before which all nodes have updated their states and proceed to perform transmissions to their neighboring nodes. 
The upper bound $\mathcal{B}$ is translated to an upper bound of $T$ in actual time (in seconds). 
This is mainly because it is not possible for nodes to count the number of time steps elapsed in the network (and understand when $\mathcal{B}$ time steps have passed. 
The value $T$ can be counted by each node individually.}


\subsection{Asynchronous $\max$/$\min$ - Consensus}\label{Prel_asy_maxmincons}

In asynchronous max/min consensus (see \cite{2013:Giannini}), the update rule for every node $v_i \in \mathcal{V}$ is: 
\begin{align}\label{asynch_max_operation_eq}
x_i^{[k + \theta_i^{[k]}]} = \max_{v_{j}\in \mathcal{N}_i^{-} \cup \{v_{i}\}}\{ x_j^{[k + \theta_{ij}^{[k]}]} \}, 
\end{align}
where $\theta_i^{[k]}$ is the update instance of node $v_i$, $x_j^{[k + \theta_{ij}^{[k]}]}$ are the states of the in-neighbors $v_{j}\in \mathcal{N}_i^{-} \cup \{v_{i}\}$ during the time instant of $v_i$'s update, $\theta_{ij}^{[k]}$ are the asynchronous state updates of the in-neighbors of node $v_i$ that occur between two consecutive updates of node $v_i$'s state. 
The asynchronous max/min consensus in \eqref{asynch_max_operation_eq} converges to the maximum value among all nodes in a finite number of steps $s' \leq D \mathcal{B}$ (see \cite{2013:Giannini}), where $D$ is the diameter of the network, and $\mathcal{B}$ is the upper bound on the number of time steps required for a node $v_j$ to process the information received from its in-neighbors.

%
%
%
%
\section{Distributed Asynchronous Optimization via ADMM with Efficient Communication}\label{sec:distr_ADMM_quant}

In this section we present a distributed algorithm which solves problem \eqref{objective_function_no_quant}. 
Before presenting the operation of the proposed algorithm, we analyze the ADMM operation over the problem \eqref{objective_function_no_quant}.  

In \eqref{objective_function_no_quant}, let us denote $F(X) \coloneqq \sum_{i=1}^{n} f_i(x_i)$. 
This means that the Lagrangian function is equal to
\begin{equation}\label{Lagrangian}
	L(X,z,\lambda) = F(X) + g(z) + \lambda\T(X - z),
	\end{equation}
where $ \lambda  \in  \mathbb{R}^{np} $ is the Lagrange multiplier. 
We now make the following assumptions to solve the problem \eqref{objective_function_no_quant}. 

\begin{assum}\label{assup_convex}
		Every cost function $f_i:\mathbb{R}^{p} \rightarrow \mathbb{R}$ is closed, proper and convex.
\end{assum}

\begin{assum}\label{assup_saddel_point}
		The Lagrangian $ L(X,z,\lambda) $ has a saddle point. 
        This means that there exists $(X^{*},z^{*},\lambda^{*})$, for which 
		\begin{equation}\label{saddle_point}
		      L(X^{*},z^{*}, \lambda)\le L (X^{*},z^{*},\lambda^{*})\le L(X,z,\lambda^{*}),
		\end{equation}
		for all $X \in \mathbb{R}^{np}$, $z \in \mathbb{R}^{np}$, and $\lambda \in \mathbb{R}^{np}$.
\end{assum}

Assumption~\ref{assup_convex} means that the local cost function $f_i$ of every node $v_i$ can be non-differentiable (see \cite{boyd2011distributed}). 
Furthermore, Assumptions~\ref{assup_convex} and~\ref{assup_saddel_point} mean that $ L(X,z,\lambda^{*}) $ is convex in $(X,z)$ and $ (X^{*},z^{*}) $ is a solution to problem~\eqref{objective_function_no_quant} (see \cite{boyd2011distributed,wei2012distributed}). 
Note that this is also based on the definition of $g(z)$ in~\eqref{indicator_function_g}. 
{Note here that our results extend naturally to strongly convex cost functions, since strong convexity implies convexity.} 


Let us now focus on the Lagrangian of the problem in \eqref{objective_function_no_quant}. 
At time step $k$, the augmented Lagrangian of \eqref{objective_function_no_quant} is 
\begin{align}\label{augmented_Lagrangian2}
	L_{\rho}&(X^{[k]},z^{[k]},\lambda^{[k]}) \\
	=& \sum_{i=1}^{n} f_i(x_i^{[k]}) + g(z^{[k]}) + {\lambda^{[k]}}\T(X^{[k]} - z^{[k]} ) \nonumber \\
        +& \frac{\rho}{2}\|X^{[k]} - z^{[k]} \|^{2}    \nonumber \\
	=&  \sum_{i=1}^{n}\left( f_i(x_i^{[k]}) + {\lambda_i^{[k]}}\T(x_i^{[k]} - z_i^{[k]}) + \frac{\rho}{2}\|x_i^{[k]} - z_i^{[k]} \|^{2} \right)  \nonumber \\
        +& g(z^{[k]}) ,   \nonumber
	\end{align}
	where $ z_i\in \mathbb{R}^{p}  $ is the $i^{th}$ element of vector $z$. 
In \eqref{admm_iterations_x}--\eqref{admm_iterations_lambda} we ignore terms that are independent of the optimization variables such as $x_i, z$ for node $v_i$. 
This means that \eqref{admm_iterations_x}--\eqref{admm_iterations_lambda} become: 
\begin{align}
	x_i^{[k+1]} =&  \operatorname*{argmin}_{x_i} f_i(x_i) + {\lambda_i^{[k]}}\T x_i + \frac{\rho}{2}\|x_i - z_i^{[k]} \|^{2}, \label{dadmm_x}\\
	z^{[k+1]} =& \operatorname*{argmin}_z g(z) + {\lambda^{[k]}}\T(X^{[k+1]} - z )  + \frac{\rho}{2}\|X^{[k+1]} - z \|^{2}\nonumber \\
	=& \operatorname*{argmin}_z g(z) +\frac{\rho}{2}\|X^{[k+1]} - z + \frac{1}{\rho}\lambda^{[k]} \|^{2}, \label{dadmm_z}\\
	\lambda_i^{[k+1]} =& \lambda_i^{[k]} + \rho (x_i^{[k+1]} - z_i^{[k+1]}) \label{dadmm_lamda}. 
	\end{align}
Note that for \eqref{dadmm_z} we use the identity $2a^Tb + b^2 = (a+b)^2-a^2$ for $a = \lambda^{[k]}/\rho$ and $b = X^{[k+1]} - z$. 

Equations \eqref{dadmm_x}, \eqref{dadmm_lamda} can be executed independently by node $v_i$ in a parallel fashion. 
Specifically, node $v_i$ can solve \eqref{dadmm_x} for $x_i^{[k+1]}$ by a classical method (e.g., a proximity operator~\cite[Section 4]{boyd2011distributed}), and implement trivially \eqref{dadmm_lamda} for $\lambda_i^{[k+1]}$. 

In \eqref{indicator_function_g}, $g(z)$ is the indicator function of the closed nonempty convex set $\mathcal{C}$.
This means that \eqref{dadmm_z} becomes 
\begin{align}
    z^{[k+1]} = \Pi_{\mathcal{C}}(X^{[k+1]}+ \lambda^{[k]}/\rho) , 
\end{align}
where $\Pi_{\mathcal{C}}$ is the projection (in the Euclidean norm) onto $\mathcal{C}$.
It is important to note here that the elements of $z$ (i.e., $ z_1, z_2, \ldots, z_n $) should belong into the set $\mathcal{C}$ in finite time. 
This is due to the definition of $g(z)$ in~\eqref{indicator_function_g}. 
Specifically, if the elements of $z$ do not belong in $\mathcal{C}$ then $g(z)  = \infty$ (thus \eqref{dadmm_z} cannot be executed). 
Therefore, we need to adjust the values of the elements of $z$ so  that they belong in the set $\mathcal{C}$ in finite time (i.e., we need to set the elements of $z$ such that $\| z_i - z_j\| \le \epsilon, \forall v_i, v_j \in \mathcal{V}$). 
Note that if $z_i - z_j = 0, \forall v_i, v_j \in \mathcal{V}$, then every node $v_i \in \mathcal{V}$ has reached consensus.
Specifically, we can have that in finite time the state $z_i$ becomes
\begin{align}\label{consensus_z}
    z_i = \frac{1}{n} \sum_{l=1}^{n}z_l^{[0]} , \forall v_i \in \mathcal{V} ,
\end{align}
where $z_l^{[0]} = x_l^{[k+1]} + \lambda_l^{[k]}/\rho$. 
Furthermore, $\| z_i - z_j\| \le \epsilon, \forall v_i, v_j\in\mathcal{V}$ means that 
\begin{align}\label{consensus_z2}
    z_i \in [ \frac{1}{n} \sum_{l=1}^{n}z_l^{[0]} - \frac{\epsilon}{2} ,  \frac{1}{n} \sum_{l=1}^{n}z_l^{[0]} + \frac{\epsilon}{2} ], \forall v_i \in \mathcal{V} ,
\end{align}
where $z_l^{[0]} = x_l^{[k+1]} + \lambda_l^{[k]}/\rho$. 
This means that for every node $v_i$, $z_i$ enters a circle with its center at $\frac{1}{n} \sum_{l=1}^{n} (x_l^{[k+1]} + \lambda_l^{[k]}/\rho)$ and its radius as $\epsilon / 2$.
Finally, from \eqref{objective_function_no_quant}, we have that each node in the network needs to communicate with its neighbors in an efficient manner. 
For this reason, we aim to allow each node $v_i$ coordinate with its neighboring nodes by exchanging quantized values in order to fulfil \eqref{consensus_z2}.

\subsection{Distributed Optimization Algorithm}\label{alg_ADMM_quant}

We now present our distributed optimization algorithm. 
The algorithm is detailed below as Algorithm~\ref{alg1} and allows each node in the network to solve the problem presented in \eqref{objective_function_no_quant}. 
The operation of the proposed algorithm is based on two parts. 
During these parts, each node $v_i$ (i) calculates $x_i^{[k+1]}$, $z_i^{[k+1]}$, $\lambda_i^{[k+1]}$ according {to \eqref{dadmm_x}--\eqref{dadmm_lamda}} (see Algorithm~\ref{alg1}), and (ii) coordinates with other nodes in a communication efficient manner in order to calculate $z_i^{[k+1]}$ that belongs in $\mathcal{C}$ in \eqref{setC} (see Algorithm~\ref{alg2}). 
Note that Algorithm~\ref{alg2} is a finite time coordination algorithm with quantized communication and is executed as a step of Algorithm~\ref{alg1}. 

Note that during Algorithm~\ref{alg1}, nodes operate in an asynchronous fashion. 
Synchronous operation requires synchronization among nodes or the existence of a global clock so that all nodes to agree on their update time. 
In our setting, asynchronous operation arises when each node (i) starts calculating $x_i^{[k+1]}$, $z_i^{[k+1]}$, $\lambda_i^{[k+1]}$ according {to \eqref{dadmm_x}--\eqref{dadmm_lamda}} in Algorithm~\ref{alg1}, and (ii) calculates $z_i^{[k+1]}$ that belongs in $\mathcal{C}$ in \eqref{setC} in Algorithm~\ref{alg2}. 
This can be achieved by making the internal clocks of all nodes have similar pacing, which will allow them to execute the optimization step at roughly the same time~\cite{lamport_time_2019}. 
Furthermore, making the internal clocks of all nodes have similar pacing does not mean that we have to synchronize the clocks of the nodes (or their time-zones). 
Note that this is a common procedure in most modern computers as the clock pacing specification is defined within the Advanced Configuration and Power Interface ({ACPI}) specification~\cite{AdvancedConfigurationandPowerInterfa}.

We now make the following assumption which is important for the operation of our algorithm. 

\begin{assum}\label{digr_diam}
The diameter $D$ (or an upper bound) is known to every node $v_i$ in the network. 
\end{assum}

Assumption~\ref{digr_diam} is necessary so that each node $v_i$ is able to determine whether calculation of $z_i$ that belongs in $\mathcal{C}$ in \eqref{setC} has been achieved in a distributed manner. 
We now present the details of Algorithm~\ref{alg1}.

\noindent
\vspace{-0.3cm}    
\begin{varalgorithm}{1}
\caption{QuAsyADMM - Quantized Asynchronous ADMM}
\textbf{Input:} Strongly connected $\mathcal{G} = (\mathcal{V}, \mathcal{E})$, parameter $\rho$, diameter $D$, error tolerance $\epsilon \in \mathbb{Q}$, upper bound on processing delays $\mathcal{B}$. 
Assumptions~\ref{assump_time_index_bound}, \ref{assup_convex}, \ref{assup_saddel_point}, \ref{digr_diam} hold. 
{$k_{\text{max}}$ (ADMM maximum number of iterations).}
\\
\textbf{Initialization:} Each node $v_i \in \mathcal{V}$ sets randomly $x^{[0]}, z^{[0]}, \lambda^{[0]}$, and sets $\Delta = \epsilon / 3$. \\ 
\textbf{Iteration:} For $k=0,1,2,\dots {, k_{\text{max}}}$, each node $v_i \in \mathcal{V}$ does the following: 
\begin{list4}
\item[1)] Calculate $x_i^{[k+1]}$ via \eqref{dadmm_x}; 
\item[2)] Calculate $z_i^{[k+1]}$ = Algorithm~\ref{alg2}($x_i^{[k+1]} + \lambda_i^{[k]}/\rho, D, \Delta, \mathcal{B}$); 
\item[3)] Calculate $\lambda_i^{[k+1]}$ via \eqref{dadmm_lamda}.  
\end{list4}
\textbf{{Output:}} Each node $v_i \in \mathcal{V}$ calculates $x_i^*$ which solves problem \eqref{objective_function_no_quant} in Section~\ref{sec:probForm}. 
\label{alg1} 
\end{varalgorithm}

\noindent
\vspace{-0.3cm}    
\begin{varalgorithm}{2}
\caption{QuAsAvCo - Quantized Asynchronous Average Consensus}
\textbf{Input:} $x_i^{[k+1]} + \lambda_i^{[k]}/\rho, D, \Delta, \mathcal{B}$. 
\\
\textbf{Initialization:} Each node $v_i \in \mathcal{V}$ does the following: 
\begin{list4}
\item[$1)$] Assigns probability $b_{li}$ to each out-neigbor $v_l \in \mathcal{N}^+_i \cup \{v_i\}$, as follows
\begin{align*}
b_{li} = \left\{ \begin{array}{ll}
         \frac{1}{1 + \mathcal{D}_i^+}, & \mbox{if $l = i$ or $v_{l} \in \mathcal{N}_i^+$,} \\
         0, & \mbox{if $l \neq i$ and $v_{l} \notin \mathcal{N}_i^+$;}\end{array} \right. 
\end{align*} 
\item[$2)$] $\text{flag}_i = 0$, $\xi_i = 2$, $y_i = 2 \  q_{\Delta}^a(x_i^{[k+1]} + \lambda_i^{[k]}/\rho)$ (see \eqref{asy_quant_defn}); 
\end{list4} 
\textbf{Iteration:} For $\eta = 1,2,\dots$, each node $v_i \in \mathcal{V}$, does: 
\begin{list4a}
\item[$1)$] \textbf{if} $\eta \mod (D \mathcal{B}) = 1$ \textbf{then} sets $M_i = \lceil y_i / \xi_i \rceil$, $m_i = \lfloor y_i / \xi_i \rfloor$; 
\item[$2)$] broadcasts $M_i$, $m_i$ to every $v_{l} \in \mathcal{N}_i^+$; receives $M_j$, $m_j$ from every $v_{j} \in \mathcal{N}_i^-$; sets $M_i = \max_{v_{j} \in \mathcal{N}_i^-\cup \{ v_i \}} M_j$, \\ $m_i = \min_{v_{j} \in \mathcal{N}_i^-\cup \{ v_i \}} m_j$; 
\item[$3)$] sets $d_i = \xi_i$; 
\item[$4)$] \textbf{while} $d_{i} > 1$ \textbf{do} 
\begin{list4a}
\item[$4.1)$] $c_i^{[\eta]} = \lfloor y_{i} \  / \  \xi_{i} \rfloor$; 
\item[$4.2)$] sets $y_{i} = y_{i} - c_i^{[\eta]}$, $\xi_{i} = \xi_{i} - 1$, and $d_i = d_i - 1$; 
\item[$4.3)$] transmits $c_i^{[\eta]}$ to randomly chosen out-neighbor $v_l \in \mathcal{N}^+_i \cup \{v_i\}$ according to $b_{li}$; 
\item[$4.4)$] receives $c_j^{[\eta]}$ from $v_j \in \mathcal{N}_i^-$ and sets 
\begin{align}
y_i & = y_i + \sum_{j=1}^{n} \sum_{r=0}^{\mathcal{B}}  w^{[r]}_{\eta-r,ij} \ c^{[\eta-r]}_{j} \ , \\
\xi_i & = \xi_i + \sum_{j=1}^{n} \sum_{r=0}^{\mathcal{B}} w^{[r]}_{\eta-r,ij} \ ,
\end{align}
where $w^{[r]}_{\eta-r,ij} = 1$ when the processing time of node $v_i$ is equal to $r$ at time step $\eta-r$, so that node $v_i$ receives $c^{[\eta]}_{i}$, $1$ from $v_j$ at time step $\eta$ (otherwise $w^{[r]}_{\eta-r,ij} = 0$ and $v_i$ receives no message at time step $\eta$ from $v_j$);
\end{list4a} 
\item[$5)$] \textbf{if} $\eta \mod (D \mathcal{B}) = 0$ \textbf{and} $M_i - m_i \leq 1$ \textbf{then} sets $z_i^{[k+1]} = m_i \Delta$ and stops operation. 
\end{list4a} 
\textbf{Output:} $z_i^{[k+1]}$.
\label{alg2} 
\end{varalgorithm}

Note here that Algorithm~\ref{alg1} shares similarities with \cite{2021:Wei_Themis}. 
However, Algorithm~\ref{alg1} is adjusted to guarantee efficient communication between nodes by allowing them to exchange quatized valued messages (this characteristic is not present in \cite{2021:Wei_Themis}). 
More specifically, the intuition of Algorithm~\ref{alg1} is the following. 
Each node $v_i$ aims to solve the problem presented in \eqref{objective_function_no_quant} via the ADMM strategy. 
Initially, each node chooses a suitable quantization level $\Delta$ so that the constraints of \eqref{objective_function_no_quant} are fulfilled (as explained later in Remark~\ref{epsilon_choice}).
During each time step $k$ each node $v_i$ calculates $x_i^{[k+1]}$ via \eqref{dadmm_x}. 
Then, each node $v_i$ executes Algorithm~\ref{alg2} in order to calculate $z_i^{[k+1]}$ which belongs in the set $\mathcal{C}$ (i.e., $z_i^{[k+1]} \in \mathcal{C}, \forall v_i \in \mathcal{V}$). 
Finally, each node $v_i$ uses $x_i^{[k+1]}$ and the result of Algorithm~\ref{alg2} in order to calculate $\lambda_i^{[k+1]}$ via \eqref{dadmm_lamda}. 
Algorithm~\ref{alg2} is a distributed coordination algorithm which allows each node to calculate the quantized average of each node's initial state. 
{The main characteristic of Algorithm~\ref{alg2} is that it combines (i) (asymmetric) quantization, (ii) quantized averaging, and (iii) a stopping strategy.  
In Algorithm~\ref{alg2}, initially each node $v_i$ uses an asymmetric quantizer to quantize its state; {see Initialization-step~$2$.} 
Then, at each time step $\eta$ each node $v_i$: 
\begin{itemize}
    \item splits the $y_i$ into $\xi_i$ equal pieces (the value of some pieces might be greater than others by one); {see Iteration-steps~$4.1$, $4.2$,}
    \item transmits each piece to a randomly selected out-neighbor or to itself; {see Iteration-step~$4.3$,}
    \item receives the pieces transmitted from its in-neighbors, sums them with $y_i$ and $\xi_i$, and repeats the operation; {see Iteration-step~$4.4$.} 
\end{itemize}
Finally, every $D \mathcal{B}$ time steps, each node $v_i$ performs in parallel a max-consensus and a min-consensus operation; {see Iteration-steps~$1$ and $2$.} 
If the results of the max-consensus and min-consensus have a difference less or equal to one (see Iteration-steps~$5$), each node $v_i$ (i) scales the solution according to the quantization level, (ii) stops the operation of Algorithm~\ref{alg2}, (iii) uses the value $z_i^{[k+1]}$ to continue the operation of Algorithm~\ref{alg1}.} 
{Note that Algorithm~\ref{alg2} converges in finite number of steps according to \cite[Theorem~$1$]{2021:Rikos_Hadj_Splitting_Autom}, since QuAsAvCo has similar structure to that in \cite{2021:Rikos_Hadj_Splitting_Autom}}.

\begin{remark}\label{epsilon_choice} 
It is important to note here that during the initialization of Algorithm~\ref{alg1}, the error tolerance $\epsilon$ is chosen to be a rational number (i.e., $\epsilon \in \mathbb{Q}$). 
This is not a limitation for the ADMM optimization process in Algorithm~\ref{alg1}. 
The real-valued $\epsilon$ can be chosen such that it can be represented as a rational value. 
Furthermore, this choice facilitates the operation of Algorithm~\ref{alg2}. 
Specifically, a rational value for $\epsilon$ facilitates the choice of a suitable quantization level $\Delta$ (since $\Delta = \epsilon/3$).  
During the execution of Algorithm~\ref{alg2} nodes quantize their states, thus an error $e_{q_1} \leq \Delta$ is imposed to every state. 
Then, Algorithm~\ref{alg2} converges to the quantized average thus, the final states of the nodes have an error $e_{q_2} \leq \Delta$. 
This means that after executing Algorithm~\ref{alg2}, we have $| z_i - z_j | \leq 2\Delta < \epsilon$, and thus we have $z_i^{[k+1]} \in \mathcal{C}$ in \eqref{setC}, $\forall v_i \in \mathcal{V}$.  
For this reason, any choice of $\Delta$ for which $\Delta < \epsilon / 2$ is suitable for the operation of our algorithm for a given error tolerance $\epsilon$. 
\end{remark}

\begin{remark} 
In practical applications, nodes do not know the value of $\mathcal{B}$. 
However, $B$ time-steps (which is its upper bound) is guaranteed to be executed within $T$ seconds (see Assumption~\ref{assump_time_index_bound}). 
As noted previously, consistent pacing of each node's clock ensures that the check for convergence at each node will happen at roughly the same time (see \cite{lamport_time_2019}).
Therefore, at every $DT$ seconds, each node checks whether Algorithm~\ref{alg2} can be terminated. 
\end{remark}

\subsection{Convergence of Algorithm~\ref{alg1}}\label{ConvADMMAlg}

We now analyze the convergence time of Algorithm~\ref{alg1} via the following theorem. 
{Our theorem is inspired from \cite{2021:Wei_Themis} but is adjusted to the quantized nature of Algorithm~\ref{alg1}. 
However, due to space limitations we omit the proof (we will include it at an extended version of our paper).}



\begin{theorem}\label{converge_Alg1}
Let us consider a strongly connected digraph $\mathcal{G} = (\mathcal{V}, \mathcal{E})$. 
Each node $v_i \in \mathcal{V}$, is endowed with a scalar local cost function $f_i(x): \mathbb{R}^p \mapsto \mathbb{R}$, and Assumptions~\ref{assump_time_index_bound}-\ref{digr_diam} hold. 
Furthermore, every node $v_i$ has knowledge of a parameter $\rho$, the network diameter $D$, an error tolerance $\epsilon \in \mathbb{Q}$, and an upper bound on processing delays $\mathcal{B}$. 
During the operation of Algorithm~\ref{alg1}, let us consider the variables $\{X^{[k]}, z^{[k]}, \lambda^{[k]}\}$, where $X^{[k]} = [{x_1^{[k]}}\T,{x_2^{[k]}}\T, \ldots, {x_n^{[k]}}\T]\T$ and $\lambda^{[k]} = [{\lambda_1^{[k]}}\T,{\lambda_2^{[k]}}\T, \ldots, {\lambda_n^{[k]}}\T]\T$;
then, define $\bar{X}^{[k]} = \frac{1}{k} \sum_{s=0}^{k-1}X^{[s+1]}, \bar{z}^{[k]} = \frac{1}{k} \sum_{s=0}^{k-1}z^{[s+1]}$. 
During the operation of Algorithm~\ref{alg1} we have 
\begin{align}
0&\le L(\bar X^{[k]},\bar z^{[k]}, \lambda^{*})- L ( X^{*},z^{*},\lambda^{*}) \label{convergence_relationship}\\
&\le \frac{1}{k}\left( \frac{1}{2\rho} \|\lambda^{*}-\lambda^{[0]}\|^2 + \frac{\rho}{2}\|X^{*}-z^{[0]}\|^2\right) + \mathcal{O}(
2\Delta \sqrt{n}) ,\nonumber
\end{align}
for every time step $k$, where $\Delta$ is the quantization level for calculating $z_i \in \mathcal{C}$ in \eqref{setC} during the operation of Algorithm~\ref{alg2}. 
\end{theorem}

It is important to note that in Theorem~\ref{converge_Alg1} we focus on the convergence of the optimization steps, i.e., the steps executed during the operation of Algorithm~\ref{alg1}. 
Due to the operation of Algorithm~\ref{alg2} we have that in \eqref{convergence_relationship} an additional term $\mathcal{O}(2\Delta \sqrt{n})$ appears. 
This term (as will be seen later in Section~\ref{sec:results}) affects the precision according to which the optimal solution is calculated. 
However, we can adjust Algorithm~\ref{alg2} to operate with a dynamically refined quantization level $\Delta$. 
For example, we can initially set $\Delta = \epsilon / 3$ (where $\epsilon \in \mathbb{Q}$). 
Then, execute Algorithm~\ref{alg2} during every time step $k$ with quantization level $\Delta' = \frac{\Delta}{10(k+1)}$. 
Since we have $\frac{\Delta}{10(k+1)} < \frac{\Delta}{10(k)}$ for every $k$, then, Algorithm~\ref{alg2} will lead to a reduction of the error on the optimal solution that depends on the quantization level (i.e., the term $\mathcal{O}(2\Delta \sqrt{n})$ in \eqref{convergence_relationship} will be reduced after every execution of Algorithm~\ref{alg2}).
However, please note that this analysis is outside of the scope of this paper and will be considered in an extended version. 

%
%
%
%

\section{Simulation Results} \label{sec:results}
    
In this section, we present simulation results in order to demonstrate the operation of Algorithm~\ref{alg1} and its advantages. 
Furthermore, we compare Algorithm~\ref{alg1} against existing algorithms and  emphasize on the introduced improvements. 



\begin{figure}[t]
\begin{center}
\includegraphics[width=.95\columnwidth]{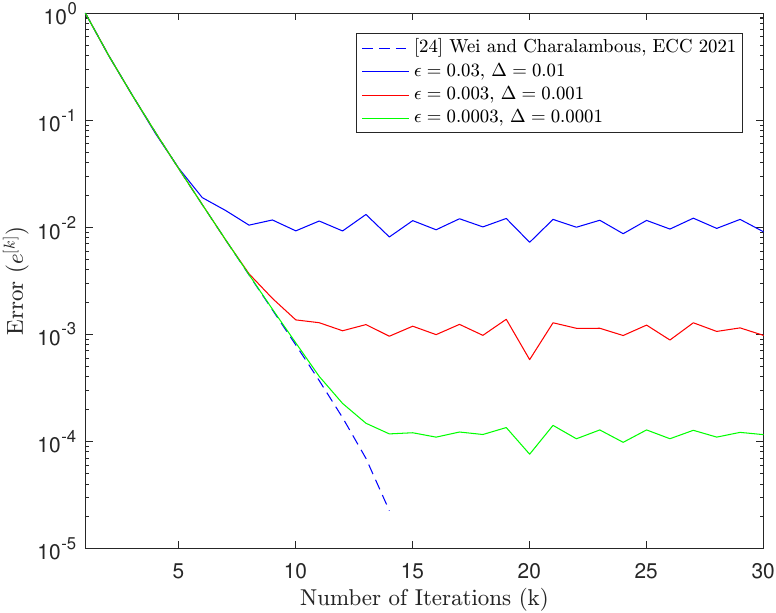}
\caption{Comparison of Algorithm~\ref{alg1} with \cite{2021:Wei_Charalambous} over a directed graph comprised of $100$ nodes for $\epsilon = 0.03, 0.003, 0.0003$, and $\Delta = 0.01, 0.001, 0.0001$, respectively\vspace{-0.2cm}.}
\label{Themis_1}
\end{center}
\end{figure}

In Fig.~\ref{Themis_1}, we focus on a network comprised of $100$ nodes modelled as a directed graph. 
{Each node $v_i$ is endowed with a scalar local cost function $f_i(x) = 0.5 x^\top P_i x + q_i^\top x + r_i$. 
This cost function is quadratic and convex. 
Furthermore, for $f_i(x)$ we have that (i) $P_i$ was initialized as the square of a randomly generated symmetric matrix $A_i$ (ensuring it is is positive definite), (ii) $q_i$ is initialized as the negation of the product of the transpose of $A_i$ and a randomly generated vector $b_i$ (i.e., it is a linear term), (iii) and $r_i$ is initialized as half of the squared norm of the randomly generated vector $b_i$ (i.e., it is a scalar constant).} 
We execute Algorithm~\ref{alg1} and we show how the nodes’ states converge to the optimal solution for $\epsilon = 0.03, 0.003, 0.0003$, and $\Delta = 0.01, 0.001, 0.0001$, respectively. 
We plot the error $e^{[k]}$ defined as 
\begin{equation}\label{eq:distance_optimal}
    e^{[k]} =  \frac{\sqrt{\sum_{j=1}^n (x_j^{[k]} - x^*)^2}}{\sqrt{\sum_{j=1}^n (x_j^{[0]} - x^*)^2}} , 
\end{equation}
where $x^*$ is the optimal solution of the optimization problem in $\eqref{objective_function_no_quant}$.
Note that from Remark~\ref{epsilon_choice}, we have that any $\Delta < \epsilon / 2$ is suitable for the operation of Algorithm~\ref{alg1} for a given $\epsilon$. 
In Fig.~\ref{Themis_1}, we execute Algorithm~\ref{alg1} for $\Delta = \epsilon / 3$. 
We can see that Algorithm~\ref{alg1} converges to the optimal solution for the three different values of $\epsilon$. 
However, Algorithm~\ref{alg1} is able to approximate the optimal solution with precision that depends on the quantization level (i.e., during Algorithm~\ref{alg1}, nodes are able to calculate a neighborhood of the optimal solution). 
Reducing the quantization level $\Delta$ allows calculation of the optimal solution with higher precision. 
Furthermore, we can see that after calculating the optimal solution our algorithm exhibits an oscillatory behavior due to quantized communication. 
This means quantized communication introduces nonlinearities to the consensus calculation which in turn affect the values of other parameters such as $x$ and $z$, and $\lambda$ (see iteration steps~$1$, $2$, $3$), and for this reason we have this oscillatory behavior. 
Finally, we can see that Algorithm~\ref{alg1} exhibits comparable performance with \cite{2021:Wei_Charalambous} (which is plotted until optimization step $14$) until the neighborhood of the optimal solution is calculated. 
However, in \cite{2021:Wei_Charalambous} nodes are able to exchange real-valued messages. 
Specifically, in \cite{2021:Wei_Charalambous} nodes are required to form the Hankel matrix and perform additional computations when the matrix loses rank.
This requires nodes to exchange the exact values of their states.
Therefore, the main advantage of Algorithm~\ref{alg1} compared to \cite{2021:Wei_Charalambous}, is that it exhibits comparable performance while guaranteeing efficient (quantized) communication among nodes. 


%
%
%
%

\section{Conclusions and Future Directions}\label{sec:conclusions}

In this paper, we presented an asynchronous distributed optimization algorithm which combines the Alternating Direction Method of Multipliers (ADMM) strategy with a finite time quantized averaging algorithm. 
We showed that our proposed algorithm is able to calculate the optimal solution while operating over directed communication networks in an asynchronous fashion, and guaranteeing efficient (quantized) communication between nodes.  
{We analyzed the operation of our algorithm and showed that it converges to a neighborhood of the optimal solution (that depends on the quantization level) at a rate of $O(1/k)$.} 
Finally, we demonstrated the operation of our algorithm and compared it against other algorithms from the literature. 

In the future, we aim to enhance the operation of our algorithm to avoid the oscillatory behavior after calculating the optimal solution. 
Furthermore, we plan to develop strategies that allow calculation of the \textit{exact} optimal solution while guaranteeing efficient communication among nodes. 
Finally, we will focus on designing efficient communication strategies for non-convex distributed optimization problems.

\bibliographystyle{IEEEtran}
\bibliography{bibliografia_consensus}

\end{document}